# Railway Virtual Coupling: A Survey of Emerging Control Techniques

Qing Wu, *Member, IEEE*, Xiaohua Ge, *Senior Member, IEEE*, Qing-Long Han, *Fellow, IEEE*, and Yafei Liu, *Member IEEE*

*Abstract*—This paper provides a systematic review of emerging control techniques used for railway Virtual Coupling (VC) studies. Train motion models are first reviewed, including model formulations and the force elements involved. Control objectives and typical design constraints are then elaborated. Next, the existing VC control techniques are surveyed and classified into five groups: consensus-based control, model prediction control, sliding mode control, machine learning-based control, and constraints-following control. Their advantages and disadvantages for VC applications are also discussed in detail. Furthermore, several future studies for achieving better controller development and implementation, respectively, are presented. The purposes of this survey are to help researchers to achieve a better systematic understanding regarding VC control, to spark more research into VC and to further speed-up the realization of this emerging technology in railway and other relevant fields such as road vehicles.

*Index Terms*— Virtual coupling; train motion model; gap references; consensus control; model prediction control; sliding mode control; machine learning

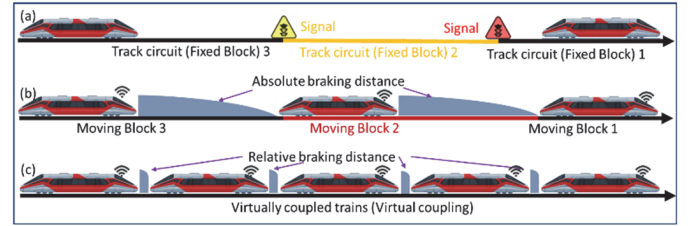

Fig. 1. Railway signaling models: (a) Fixed Block, (b) Moving Block, and (c) Virtual Coupling

## I. INTRODUCTION

AN important objective of railway signaling is to keep trains that are or will be running on the same section of track separated at a safe distance. Many of the current railway signaling systems are using the Fixed Block Signaling (FBS) model as shown in Fig. 1(a). This model is provably reliable from today's technological perspective. However, it is obviously not efficient for railway traffic as a significant amount of track space is unoccupied between adjacent running trains. With the ever-increasing demands from passenger and freight transport, enabled by technological advances in smart sensors and wireless communications, the Moving Block Signaling (MBS) model as shown in Fig. 1(b) has been developed and implemented on some railways. The MBS model uses on-board signaling systems rather than the way-side version of the FBS model. Movement authority blocks move with the trains rather than being fixed on the track. In this way, the distance between any two adjacent running trains can be significantly decreased. The shortest distance can be the absolute braking distance (ABD) of the following train plus a certain safety margin. One can easily appreciate that traffic efficiency can be significantly improved by implementing the MBS model.

Over the past several years, an exciting railway signaling model called Virtual Coupling (VC), as shown in Fig. 1(c), has attracted enormous interest from railway industry and academia. It was proposed at the end of the 1990s [1,2]. The idea is to run trains in the MBS model with a relative braking distance (RBD) rather than the ABD. Earlier research was carried out by researchers from Technical University Braunschweig [3,4] and University of Paderborn [5,6]. VC studies were not very active before 2015 partially due to the limitations from Information and Communication Technologies [7-9]. With recent advances in those areas and the driving forces from the pursuit of higher capability, flexibility, and modularity, the European initiative Shift2Rail included Virtual Coupling research into its strategic master plan in 2015, which sparked a lot of research in this area [10-14]. Xun *et al.* [15] provided a list of projects that were recently funded all around the world. VC can significantly reduce train operation headways [16] to increase line capacity [17]. VC can also help to reduce energy costs under specific specialized cases [18]. Market potential studies [19] have indicated that VC train operations can be very attractive to customers of various rail transport sectors including high-speed, main-line and regional with benefits that are especially relevant for freight trains. System feasibility for VC has also been discussed and proven at various stages [20-22]. The first implementation and tests on low-speed trams were reported in 2019 [23].

The control system of connected trains can be regarded as a networked control system, which usually consists of four parts: sensors, communications, controllers, and actuators. This paper focuses on the controller part which has also been identified as a critical step towards the successful implementation of VC operations [21,22,24]. Autonomous driving and control for intelligent road vehicles [25-29] have been well studied and can contribute to the control development for railway VC. However, railway trains have larger and more complex systems

Corresponding author: Qing-Long Han
Q. Wu is with Centre for Railway Engineering, Central Queensland University, Rockhampton, QLD 4701, Australia (e-mail: q.wu@cqu.edu.au).

X. Ge and Q.-L. Han are with School of Science, Computing and Engineering Technologies, Swinburne University of Technology, Melbourne, VIC 3122, Australia (e-mail: xge@swin.edu.au; qhan@swin.edu.au).
Y. Liu is with School of Transportation and Logistics, Southwest Jiaotong University, Chengdu, 610031, China (e-mail: yafei.liu@swjtu.edu.cn).

than road vehicles. Longitudinal Train Dynamics (LTD) [30] are also significantly different from longitudinal road vehicle dynamics [31] for controller design and implementation. This motivates us to initiate a systematic review regarding the emerging control techniques used for railway VC studies.

The rest of the paper is arranged as follows. Section II reviews train motion models, which lay a foundation in the VC controller development. Section III outlines the objectives, which should be observed while developing the controllers. Section IV elaborates the emerging control techniques, that have been used for VC studies. Section V discusses the advantages and disadvantages of the reviewed VC control techniques. Section VI presents some topics, that can be interesting for future VC research. Section VII concludes the paper.

## II. TRAIN MOTION MODELS

Train motion models are a fundamental part of the VC controller development. Such models can be regarded as special versions of LTD models [42]. Both train motion and LTD models have a focus on the longitudinal motions by neglecting vehicle lateral and vertical motions. A difference is that train motion models often neglect relative motions between adjacent vehicles in the same train whilst the relative motions are often an important part of LTD studies. This section reviews train motion model formulations that are used in VC controllers as well as the force elements considered in these controllers.

### A. Model Formulations

Consider a group of $N (\geq 2)$ automated trains whose longitudinal motions can be regulated in a distributed and cooperative manner. In this case, each train can be regarded as an 'agent' in the context of a multi-agent system. From an algebraic graph perspective, the communication topology can be modelled by a generic digraph $\mathcal{G} = \{\mathcal{V}, \mathcal{A}, \mathcal{E}\}$, where $\mathcal{V} = \{1, 2, \cdots, N\}$ is the node (or vertex) set, $\mathcal{A} = [a_{ij}]_{N \times N}$ is the adjacency matrix with nonnegative adjacency elements (or called coupling gains) $a_{ij}$, and $\mathcal{E} \subseteq \mathcal{V} \times \mathcal{V}$ is the edge set of paired nodes with $(j, i)$ denoting an edge rooted at node $j$ and ended at node $i$. For any $i, j \in \mathcal{V}$, the adjacency element $a_{ij} > 0$ means that there is an information link from node $j$ to node $i$; and $a_{ij} = 0$, otherwise. Reaching this step, a typical motion model for each train $i \in \mathcal{V}$ can be given as [32,33]:

$$\dot{x}_{i,t} = v_{i,t} \tag{2.1}$$
$$(1 + \gamma)M_i \dot{v}_{i,t} = F_{i,t}^{db} - F_{i,t}^{re1} - F_{i,t}^{re2} - F_{i,t}^{re3} - \omega_{i,t} \tag{2.2}$$
$$\begin{cases} F_{i,t}^{re1} = (k_{1,i} + k_{2,i}v_{i,t} + k_{3,i}k_{4i}v_{i,t}^2)M_i g \\ F_{i,t}^{re2} = (k_{5i}/R_c(x_{i,t}))M_i g \\ F_{i,t}^{re3} = \sin(\theta(x_{i,t}))M_i g \end{cases} \tag{2.3}$$

where $x_{i,t}$ and $v_{i,t}$ denote the longitudinal position and velocity of train $i$ at the continuous time $t \in \mathbb{R}_{\geq 0}$, respectively; $\gamma$ represents the effect of rotational inertia of rotational components such as wheelsets and motor rotors; $M_i = M_0 + \widetilde{M}_i$ denotes the unknown mass of the train with $M_0$ being the nominal (measured) part and $\widetilde{M}_i$ being the uncertain part of the mass, respectively; $F_{i,t}^{db}$ is the tractive/brake (T/B) force; $F_{i,t}^{re1}$ is the rolling resistance; $F_{i,t}^{re2}$ is the curving resistance; $F_{i,t}^{re3}$ is the track gradient force; $\omega_{i,t}$ stands for the other uncertain inputs that were not specifically modelled by the previous components; $k_{1,i}, k_{2,i}, k_{3,i}$ are the basic empirical parameters for train rolling resistance; $k_{4,i}$ is an extra resistance parameter for tunnel resistance; $g$ is the gravity constant; $k_{5i}$ is curving resistance coefficient; $R_c(x_{i,t})$ is track curve radius; and $\theta(x_{i,t})$ is the track gradient. Note that, when the track profile can be identified beforehand and the real-time train speed information is available, one may further define a normalized control (acceleration) input $u_{i,t} = (F_{i,t}^{db} - F_{i,t}^{re1} - F_{i,t}^{re2} - F_{i,t}^{re3})/((1 + \gamma)M_0)$ and a lumped unknown input $w_{i,t} = -((1 + \gamma)\widetilde{M}_i \dot{v}_{i,t} + \omega_{i,t})/((1 + \gamma)M_0)$. Then, the motion model above can be formulated as the following linear second-order state-space model:

$$\begin{aligned} \dot{x}_{i,t} = v_{i,t}, \; \dot{v}_{i,t} = u_{i,t} + w_{i,t}, \; i \in \mathcal{V} \quad \Leftrightarrow \\ \dot{s}_{i,t} = As_{i,t} + B(u_{i,t} + w_{i,t}), \; i \in \mathcal{V} \end{aligned} \tag{2.4}$$

where $s_{i,t} = [x_{i,t}, v_{i,t}]^T, A = [0,1; 0,0]$, and $B = [0; 1]$. The primary control objective is then to design a high-level cooperative longitudinal control law $u_{i,t}$ for each automated train $i \in \mathcal{V}$. Once $u_{i,t}$ is determined, the actual low-level T/B force $F_{i,t}^{db}$ can be calculated based on the pre-saved track information and resistance characteristics.

The train motion model above is a second-order model which includes information of train position and velocity only. Table I lists the train motion models used for various VC studies. Most existing studies used the second-order motion models. Alternatively, there are also third-order state-space train motion models that can be used to incorporate more information into the train motion models. In addition to (2.1) and (2.2), some third-order models also use the following force dynamics [34]:

$$\tau_i \dot{F}_{i,t}^{db} = u_{i,t} - F_{i,t}^{db} \tag{2.5}$$

where $\tau_i = \tau_0 + \tilde{\tau}_i$ denotes the uncertain inertial lag (in unit of seconds) for the train motions, and $u_{i,t}$ represents the actual T/B input, which is regarded as the desired control input. If one further denotes $s_{i,t} = [x_{i,t}, v_{i,t}, F_{i,t}^{db}]^T, A = [0,1,0; 0,0,1/((1 + \gamma)M_0); 0,0, -1/\tau_0], B = [0; 0; 1/\tau_0], E = [0; -1/(1 + \gamma)M_0; 0], F = [0,0; -1/((1 + \gamma)M_0), 0; 0, -1/\tau_0], f_{i,t} = F_{i,t}^{re1} + F_{i,t}^{re2} + F_{i,t}^{re3}$ and $w_{i,t} = [(1 + \gamma)\widetilde{M}_i \dot{v}_{i,t} + \omega_{i,t}, \tilde{\tau}_i \dot{F}_{i,t}^{db}]^T$, the following nonlinear third-order state-space model can be derived:

$$\dot{s}_{i,t} = As_{i,t} + Bu_{i,t} + Ef_{i,t} + Fw_{i,t}, i \in \mathcal{V} \tag{2.6}$$

where the system matrices $A, B, E, F$ are known constants but the nonlinear resistances $f_{i,t}$ and the lumped uncertainties $w_{i,t}$ allow to be generally unknown. Naturally, the existence of the unknown and nonlinear inputs $f_{i,t}$ and $w_{i,t}$ makes the above train motion model more realistic and comprehensive. However, they also pose a significant challenge to controller design. In general, some dedicated control strategies, such as nonlinear train control and adaptive train control, are required to deal with the analysis and synthesis challenges.

TABLE I TRAIN MOTION MODELS FOR VC STUDIES (T/B: traction and brake force; RI: rotational inertia; RR: rolling resistance; TR: tunnel resistance; CR: curve resistance; GF: grade force; UF: uncertain force; Y: considered; SMC: sliding mode control; MPC: model predicative control; ML: machine learning; CBC: consensus-based control; and CFC: constraint following control)

| Ref. | Order | Force elements | | | | | | | Controller |
|---|---|---|---|---|---|---|---|---|---|
| | | T/B | RR | RI | TR | CR | GF | UF | |
| [45] | 2nd | | | | Y | Y | Y | Y | SMC |
| [46] | 2nd | | | | | | | Y | SMC |
| [47] | 2nd | | | | Y | Y | Y | Y | SMC |
| [44, 47-50] | 2nd | | | | Y | Y | Y | | MPC |
| [43] | 3rd | Yes to all cells of this column | Yes to all cells of this column | | Y | Y | Y | | MPC |
| [38, 51-54] | 2nd | | | | | | Y | | MPC |
| [34] | 3rd | | | | | | | | MPC |
| [55] | 2nd | | | | | Y | Y | Y | MPC |
| [56] | 2nd | | | | | | | | MPC |
| [57] | 2nd | | | | | Y | Y | | MPC |
| [39] | 2nd | | | | | | Y | | ML |
| [37] | 2nd | | | | | Y | Y | | ML |
| [32, 33] | 2nd | | | Y | | Y | Y | | CBC |
| [58-60] | 2nd | | | | | | | | CBC |
| [64] | 2nd | | | | | | | Y | CBC |
| [61] | 3rd | | | | | | | | CBC |
| [36] | 3rd | | | | | Y | Y | | CBC |
| [62, 63] | 2nd | | | Y | Y | Y | | | CFC |

For convenience of analysis and design, the force term $\dot{F}_{i,t}^{db}$ can be further replaced by acceleration via taking the time-derivative on both sides of (2.2) and combining (2.5) and then applying the exact feedback linearization technique [35]. One then obtains the following linear third-order 'position-velocity-acceleration' train motion model [36]:

$$\dot{x}_{i,t} = v_{i,t} , \dot{v}_{i,t} = a_{i,t}, \tau_i \dot{a}_{i,t} = u_{i,t} - a_{i,t} - \omega_{i,t} \quad (2.7)$$

where $u_{i,t}$ represents the desired acceleration command to be designed. Similarly, a compact third-order state-space model can be expressed as

$$\dot{s}_{i,t} = A s_{i,t} + B(u_{i,t} + w_{i,t}), i \in \mathcal{V} \quad (2.8)$$

where $s_{i,t} = [x_{i,t}, v_{i,t}, a_{i,t}]^T, A = [0,1,0; 0,0,1; 0,0,-1/\tau_0]$, $B = [0; 0; 1/\tau_0]$, $w_{i,t} = -\tilde{\tau}_i \dot{a}_{i,t} - \omega_{i,t}$.

Comparatively, the nonlinear state-space motion models (2.6) and linear state-space motion models (2.4)/(2.8) offer a trade-off between the accuracy of train motion modelling and the simplicity of train motion controller design. Usually, the nonlinear models have better accuracy for train motion modelling; but they also rely on the proper and accurate handling (e.g., estimation and/or compensation) of the unknown and nonlinear inputs. On the other hand, linear motion models generally facilitate the analysis and synthesis of the train control systems. However, they necessitate the accuracy of *a priori* knowledge regarding the track information and train system parameters. If such information is not available or not accurate, such motion models become less instructional or even inapplicable.

Fig. 2. Considerations for train motion models: (a) tractive force limits, and (b) influences of rotational inertia (simulated train speed during emergency brake) [41]

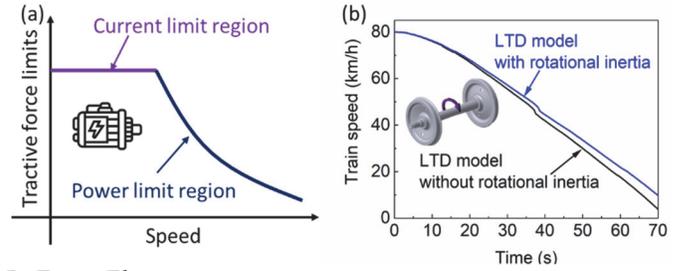

*B. Force Elements*

From the LTD perspective, the force elements considered include in-train forces, tractive forces, dynamic brake forces, air brake forces, rolling resistance (can include tunnel resistance), curving resistance and gradient forces. As discussed before, due to the neglecting of inter-vehicle motions in the VC controller development, train motion models usually do not consider in-train forces. Meanwhile, train controllers generally lump tractive forces, dynamic braking forces and air brake forces together as a single T/B force.

Table I lists the force components that were considered in different VC studies. All the models considered T/B forces and rolling resistance. Most of the models had constraints for T/B forces and the constraints had different forms. For example, Su *et al.* [37] directly limited the maximum tractive forces whilst Di Meo [36] and Liu *et al.* [38] limited the allowable train acceleration. Two different limits could be used for traction and brake cases as used in [39]. Studies that did not use specific maximum force constraints during controller development were also found. It is worth mentioning that the tractive force limits have at least two different regions as shown in Fig. 2(a). The first region is limited by the maximum allowable electrical current in the traction system whilst the second region is limited by the maximum power of the system. Similar speed-dependent maximum force limits also exist for dynamic brake forces and frictional brake forces as discussed in [40-42]. Felez *et al.* [43] used two limits for the VC controller: one for the maximum forces and the other for the maximum power (the product of speed and force). The speed-dependent limit was acknowledged in [37,44]. How the limit was implemented was not clear from the publications.

Several studies have considered tunnel resistance by adding an extra term into the rolling resistance formula. Most of the models have considered curving resistance and gradient forces which are important for realistic models. For some high-speed train studies, such as [61], it was assumed that the track grades and curvatures were small and negligible. A small number of models have used an uncertain force component to indirectly model forces such as tunnel resistance, curving resistance and gradient forces. Very few models have considered the influence of rotational inertia. Due to the simplifications applied in train motion models, rotational kinetic energy in components such as wheelsets and motor rotors cannot be simulated by modelling translations. However, during acceleration or deceleration processes, a portion of tractive or brake energy will be used to change the rotational motions. In train motion modelling, the influence of rotational inertia can be simulated by added an equivalent mass to the total translational mass of the system:

$$M_\gamma = k_\gamma J_\gamma / r^2 \quad (2.9)$$

where $k_\gamma$ is a relation parameter between the circumferential speed of the rotational component and the translational speed of the vehicle (for wheelset $k_\gamma = 1$); $M_\gamma$ is the equivalent rotational mass, $J_\gamma$ is the rotational inertia of the rotational component; and $r$ is the radius of the rotational component. The influence of rotational inertia on simulated train speed can be evident as shown in Fig. 2(b). Hence, it is recommended to be included in train motion models.

## III. CONTROL OBJECTIVES

Control objectives need to be defined and formulated prior the design of VC controllers. Relative train motion controls are discussed first in this section as the primary objective of the VC controller development. Various constraints that were used during the development of VC controllers are also discussed.

### A. Relative Motion Control

This section first discusses various gap references that can be used for VC operations. Then, motion control objectives under two different operational architectures (multi-agent and predecessor-follower) are discussed.

### Gap References

In a VC system, the design of the minimum gap reference needs to consider various factors such as: train location errors, speed errors, actuator lag, information latency (delay), train speeds, and train accelerations. Researchers in [37,57] introduced a method to determine a minimum reference gap by adapting the approach specified in the IEEE 1474.1 as shown in Fig. 3. Su et al. [51] and Xun et al. [65] used the following expression to determine the desired dynamic gap reference:

$$d_{i,t}^{de} = D_{safe} + \frac{v_{i,t}^2}{2a_b^-} - \frac{v_{i-1,t}^2}{2a_b^-}, \quad \forall i \in \mathcal{V} \quad (3.1)$$

where $d_{i,t}^{de}$ Is the (minimum) dynamic gap reference; $D_{safe}$ is the safety margin; and $a_b^-$ is the average brake acceleration. It considers the relative speed between two adjacent trains; and assumes that the two trains have the same brake accelerations. The brake accelerations can be changed when two trains have different brake performance or are located on different track grades. Liu et al. [52] replaced the average brake accelerations with maximum accelerations, which represented the best-case safety scenario, and the final minimum gap was shorter. Quaglietta et al. [66] proposed a design of a dynamic gap reference to consider various factors for: 1) train positioning errors, 2) communication delays, 3) train control delays, 4) environmental influences, and 5) brake performance differences. Usually, factors 1)-4) were lumped together as a safety margin. In this sense, this model is similar to (3.1) but more detailed in terms of being able to incorporate stochastic models for the factors considered.

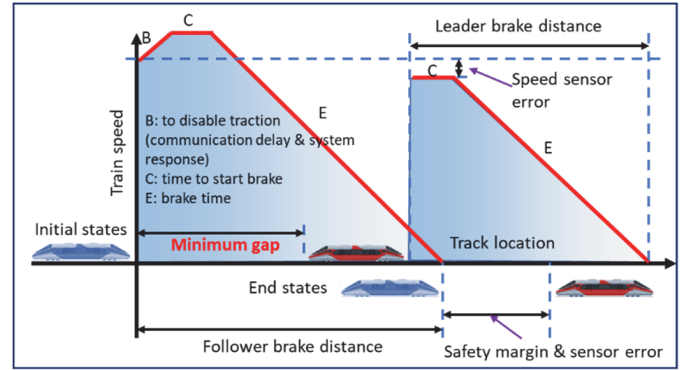

Fig. 3. Determination of the minimum gap reference

To accommodate different stages of VC operations, the gap reference also needs to achieve flexible train control operations such as train merging and splitting. Variable or transitional gap references are then needed. In [63], a transitional gap reference was designed as:

$$d_{i,t}^{dy} = k_1 \exp(-k_2 t) + d_{i,t}^{de} \quad (3.2)$$

where $d_{i,t}^{dy}$ is the transitional gap reference, $k_1$ and $k_2$ are two parameters used to adjust the transitional process. And the corresponding follower acceleration was expressed as (assuming leader has zero acceleration):

$$a_{i,t}^{dy} = k_3 k_4 \exp(-k_4 t) \quad (3.3)$$

where $a_{i,t}^{dy}$ is the follower acceleration; $k_3$ and $k_4$ are two parameters used to adjust the transitional process. The acceleration profile is shown as Profile 1 in Fig. 4. The controller started with a non-zero acceleration, which could be an issue from the perspective of passenger ride comfort. In [62], a sinusoidal function was used as the transitional gap reference:

$$d_{i,t}^{dy} = k_1 \sin[k_2(x_{i,t} + k_3) + k_4] \quad (3.4)$$

$$k_1 = \frac{d_{i,0} - d_{i,t}^{de}}{2}, k_2 = \frac{\pi}{x_{de} - x_{i,0}}, k_3 = \frac{x_{de} - 3x_{i,0}}{2}, k_4 = \frac{d_{i,0} + d_{i,t}^{de}}{2} \quad (3.5)$$

where $x_{i,t}$ is the position of the follower; $k_1$ to $k_4$ are parameters used to adjust the transitional process; $d_{i,0}$ is the initial gap; $x_{i,0}$ is the location where the transition process started; and $x_{de}$ is the location where the transition process was finished. This transition gap reference was indexed by follower position $x_{i,t}$; and the shape of the gap reference was a shifted half cycle of a sinusoidal function. The acceleration profile in the time domain is shown as Profile 2 in Fig. 4. The expression indicates that, during the transition process, there is a continuous change of acceleration between positive and negative values. This may require the T/B forces to change between negative to positive without meaningful transitions via zero T/B force states, i.e., idling states. From the train dynamics perspective, transition states, for example, staying at an idling state for five seconds between switches from traction to brake or from brake to traction, are required to minimize longitudinal impacts.

In [46], to ensure the completion of merge/separation before a given location while respecting the constraints on the jerk and acceleration of the trains, a time-based piece-wise speed profile was designed. An example section was expressed as:

$$v_{i,t} = v_{i,0} + 0.5 J_{max}(t_1)^2 + a_{max}(t - t_1) \quad (3.6)$$

where $v_{i,0}$ is the initial speed at the start of the transition; $J_{max}$ is the maximum jerk limit; $a_{max}$ is the maximum acceleration; and $t_1$ is the time required to reach maximum acceleration. The acceleration profile is shown as Profile 3 in Fig. 4. In this study, a meaningful idling state (about 25 seconds) was placed between changes between positive and negative acceleration. The utilization of idling state could evidently improve passenger comfort.

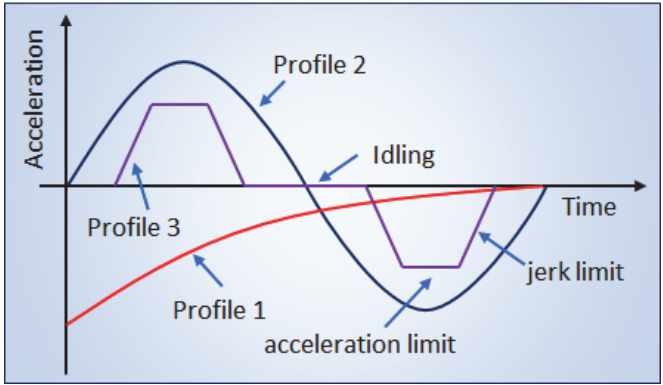

Fig. 4. Acceleration profiles for transitional gap references

*Relative Motion Control*

From the control perspective, the VC controller design aims to determine a cooperative longitudinal control law for the desired force input $F_{i,t}^{db}$ in system (2.1)-(2.2) subject to (2.3), or for the desired acceleration input $u_{i,t}$ in system (2.4) such that each follower train $i \in \mathcal{V}$ can eventually reach and maintain the same speed as the preceding train $i-1$ while preserving the desired gap reference $d_{i,t}^{de}$. Note that, for the first train $i = 1$ in the platoon, its preceding train $i = 0$ is commonly known as the train leader, which also can be modelled as

$$\dot{x}_{0,t} = v_{0,t}, (1+\gamma)M_0 \dot{v}_{0,t} = F_{0,t}^{db} + w_{0,t} \text{ or }$$
$$\dot{s}_{0,t} = A s_{0,t} + B(u_{0,t} + w_{0,t}). \quad (3.7)$$

The leader input $F_{0,t}^{db}$ or $u_{0,t}$ is usually regarded as a command reference for the entire train platoon and is usually unknown to the followers. Whether or not a follower train (except the first follower train) in the platoon can access information from the leader train is determined by the Train to Train (T2T) communication [67] topology $\mathcal{G}$. A considerable number of communication topologies used for road vehicle platooning [68-71] can also be used for VC operations. Fig. 5 shows three examples adapted from road vehicle platooning topologies.

For any $i \in \mathcal{V}$, the real-time gap between the head of train $i$ and the tail of its predecessor can be calculated as

$$d_{i,t} = x_{i,t} - x_{i-1,t} + L_{i-1}, \quad \forall i \in \mathcal{V} \quad (3.8)$$

where $L_{i-1}$ is the length of train $i$. Furthermore, the spacing error and speed error between any two consecutive trains can be described by

$$e_{i,t}^{sp} = d_{i,t} - d_{i,t}^{de}, \quad e_{i,t}^{sv} = v_{i,t} - v_{i-1,t}, \quad \forall i \in \mathcal{V}. \quad (3.9)$$

Then, the position and speed tracking errors from the train leader can be respectively formulated as:

$$e_{i,t}^p = x_{i,t} - x_{0,t} - D_{i,t}, \quad e_{i,t}^v = v_{i,t} - v_{0,t}, \quad \forall i \in \mathcal{V} \quad (3.10)$$

where $D_{i,t} = \sum_{m=1}^{i}(d_{m,t}^{de} - L_{m-1})$.

For a communication topology in which multiple trains need to be controlled at the same time, the VC controller design procedure thus becomes to determine the cooperative longitudinal control law $F_{i,t}^{db}$ or $u_{i,t}$ such that the stability of the errors $e_{i,t}^S = [e_{i,t}^{sp}, e_{i,t}^{sv}]^T$ and $e_{i,t} = [e_{i,t}^p, e_{i,t}^v]^T$ can be ensured, namely, $e_{i,t}^S \to 0$ and $e_{i,t} \to 0$ as $t \to \infty$, while also preserving driving safety, comfort and/or energy efficiency requirements.

When only two automated trains are emphasized, e.g., the follower train $i = F$ and the predecessor train $i - 1 = P$, the train motion model and controller are often known as the predecessor-follower model and predecessor-follower controller. The train motions of both the follower and predecessor can still be described by (2.1)-(2.4), while the focus of controller design will be placed on the two marshalling trains. In this case, the spacing error and speed error can be modified as:

$$e_{F,t}^{sp} = d_{F,t} - d_{F,t}^{de} = x_{F,t} - x_{P,t} + L_P - d_{F,t}^{de},$$
$$e_{F,t}^{sv} = v_{F,t} - v_{P,t}. \quad (3.11)$$

The predecessor-follower train motion model can be regarded as a special case of the multi-agent train motion model. However, due to the involvement of only two marshalling trains, the predecessor-follower model is not suitable for VC train platoons over other T2T communication topologies. In contrast, one of the prominent features of the multi-agent model is that it allows the desired VC control law to incorporate various T2T information flow structures into the controller design. In other words, the VC controller designed for the multi-agent model enjoys the benefit of learning from different neighboring trains, specified by the T2T communication topology, to regulate its motion. Furthermore, the multi-agent model empowers some further analyses of the effects of different T2T communication topologies on the train platoon stability, scalability, and control performance.

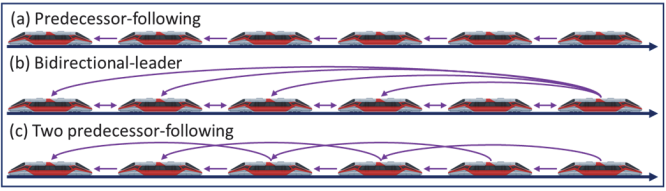

Fig. 5. VC communication topology examples: (a) Predecessor-Following, (b) Bidirectional-Leader, and (c) Two Predecessor-Following

TABLE II VC CONTROLLER CONSTRAINTS

| Ref. | Min. Gap | Max. Speed | Max. T/B (acc) | Max. Power | Max. Jerk |
|---|---|---|---|---|---|
| [43] |  | Y | Y | Y |  |
| [53,56,58] |  | Y | Y |  |  |
| [50,57] | Y | Y | Y |  | Y |
| [34, 37, 38,51,52, 54, 55,65] | Y | Y | Y |  |  |
| [39, 46, 47,59] |  |  | Y |  |  |

### B. Constraints

Various constraints were observed during the development of VC controllers as listed in TABLE II. Maximum T/B force or acceleration limits were the most used constraints. Maximum speeds and minimum following distance were also used very often. Felez *et al.* [43] applied a constraint for maximum power (the product of speed and force) of the system. Yang and Yan [57] used a maximum jerk constraint for the consideration of ride comfort. Researchers in [32,50,56] considered a stability constraint which guaranteed the actual distance between any two trains asymptotically approaches the reference; and the relative speed error between any two trains asymptotically approaches zero:

$$\lim_{t\to\infty} e_{i,t}^{sp} = \lim_{t\to\infty}(d_{i,t} - d_{i,t}^{de}) = \delta_m^+,$$
$$\lim_{t\to\infty} e_{i,t}^{sv} = \lim_{t\to\infty}(v_{i,t} - v_{i-1,t}) = \delta_m^+, \ \forall\ i \in \mathcal{V} \quad (3.12)$$

where $\delta_m^+$ is a mall value to relax the numerical pressure during the simulations and acts as zero. In addition to the references listed in TABLE II. Controllers that were not constrained were also found, e.g., [36,45,49,59].

## IV. EMERGING VIRTUAL COUPLING CONTROL TECHNIQUES

Five types of VC control techniques as listed in TABLE I are reviewed in this section. It is noted that precise classification of different control techniques is difficult. The classification breakdown used in this paper is based on the authors' research experience and also limited in number for the convenience of paper discussions. Xun *et al.* [15] also provided a classification of VC controllers, including train-following control, feedback control, optimal control, and computational intelligent control. Our controller review and classification are more detailed.

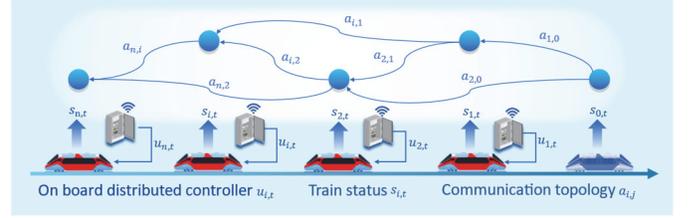

Fig. 6. Consensus-based control for VC

### A. Consensus-Based Control (CBC)

The rapid development of wireless communication technologies has greatly facilitated modern railways. Trains can communicate with each other under fast and reliable information networks. When it comes to the controller design, this means that each train can employ not only information from its direct predecessors (and/or successors) but also information from its underlying neighbors specified by a suitable T2T communication topology as shown in Fig. 6. For example, for the train motion model expressed by (2.8), a representative T2T-based VC controller for each train $i \in \mathcal{V}$ can be constructed as [32,33]:

$$u_{i,t} = K\{\sum_{j\in\mathcal{V}} a_{ij}(s_{i,t} - \tilde{d}_{i,t} - (s_{j,t} - \tilde{d}_{j,t})) + a_{i0}(s_{i,t} - \tilde{d}_{i,t} - s_{0,t})\}$$
$$= \underbrace{K\sum_{j\in\mathcal{V}\cup\{0\}} a_{ij}(s_{i,t} - \tilde{d}_{i,t})}_{\text{from train }i} \underbrace{- K\sum_{j\in\mathcal{V}} a_{ij}(s_{j,t} - \tilde{d}_{j,t})}_{\text{from neighbors of train }i} \underbrace{- a_{i0}s_{0,t}}_{\text{from leader}} \quad (4.1)$$

where $K$ denotes the controller gain matrix and $\tilde{d}_{i,t} = [D_{i,t}, 0, 0]^\top$ represents the designated gap vector. Each train $i$ can collect information from its neighboring trains specified by the nonzero adjacency elements $a_{ij}$ for any $j \in \mathcal{V}$, thus overcoming the limitation of having only information access to its direct predecessor $i-1$ in the predecessor-follower control scenario. Whether or not train $i$ can receive information from the leader train is decided by the nonzero element $a_{i0}$. Substituting (4.1) into (2.8) and combining (3.7) yields:

$$\dot{e}_t = (I_N \otimes A + \mathcal{H} \otimes BK)e_t + (I_N \otimes B)(w_t - \tilde{u}_{0,t} - \widetilde{w}_{0,t}) \quad (4.2)$$

where $e_t = [e_{1,t}, e_{2,t}, \cdots, e_{N,t}]^\top$, $w_t = [w_{1,t}, w_{2,t}, \cdots, w_{N,t}]^\top$, $\tilde{u}_{0,t} = [u_{0,t}, u_{0,t}, \cdots, u_{0,t}]^\top$, $\widetilde{w}_{0,t} = [w_{0,t}, w_{0,t}, \cdots, w_{0,t}]^\top$, $\mathcal{H} = \mathcal{L} + \mathcal{A}_0$ with $\mathcal{L}$ being the Laplacian of the graph $\mathcal{G}$ and $\mathcal{A}_0 = diag\{a_{1,0}, a_{2,0}, \cdots, a_{N,0}\}$ being the leader adjacency matrix, and $\otimes$ stands for the Kronecker product of matrices. It is obvious that the resulting closed-loop train platoon system (4.2) depends on the information of the train motion model (2.8) (characterized by $A, B$), the desired controller (specified by the gain matrix $K$), and the T2T communication topology (denoted by the matrix $\mathcal{H}$) as well as the uncertainties (denoted by $w_t, \tilde{u}_{0,t}, \widetilde{w}_{0,t}$). In this sense, the VC controller design problem can be pursued by incorporating the effects of different T2T topologies on stability, scalability, robustness, and efficiency of the train platoon system.

The above VC controller (4.1) adopts the so-called disagreement error term (the relative position, speed, acceleration differences) between train $i$ and train $j$ as feedback to regulate each train's motion. Eventually, all controlled trains will reach a consensus on the desired zero tracking errors, which implies the feasibility of the coordinated virtual coupling control objective. Such a consensus-based control strategy is motivated from the existing distributed coordinated control studies for general multi-agent systems [72] and automated vehicle platoons [69-71,73,74].

When data transmission delays $\tau_{ij,t}$ occur during T2T communications, the following distributed consensus-based coupling control protocol was proposed [36]:

$$u_{i,t} = \frac{1}{\sum_{j\in\mathcal{V}} a_{ij}} \sum_{j=0}^{N} k_{1,ij}\, a_{ij} \left[ x_{i,t} - x_{j,t-\tau_{ij,t}} - h_{ij}v_{0,t} - d_{i0}^{de} \right]$$
$$- k_{2,i}[v_{i,t} - v_{0,t}] - k_{3,i}\dot{v}_{i,t} + \frac{1}{\sum_{j\in\mathcal{V}} a_{ij}} \sum_{j=0}^{N} k_{1,ij} a_{ij}\tau_{ij,t} v_{0,t} \quad (4.3)$$

where $k_{1,ij}$, $k_{2,i}$ and $k_{3,i}$ are controller gains; $\tau_{ij,t}$ is the time-varying T2T communication delay; $h_{ij}$ is the time headway; $v_{0,t}$ is the reference speed; and $d_{i0}^{de}$ is the desired gap.

The coordinated control problem for multiple high-speed trains in the presence of actuator saturation was addressed in [59]. Specifically, the following bounded consensus-based control law was constructed to tackle the saturated input constraint:

$$u_{i,t} = \sum_{j\in\mathcal{V}} k_1 \text{sat}_\Delta(a_{ij}(e^L_{j,t} - e^L_{i,t})) - k_2 \text{sat}_\Delta(e^L_{i,t}) + k_3 \text{sat}_\Delta(f^{re1}_{i,t}) + \gamma_{i,t} \quad (4.4)$$

where $\text{sat}_\Delta(x)$ stands for the saturation function satisfying $\text{sat}_\Delta(x) = \text{sgn}(x) \cdot \Delta$ for any $|x| > \Delta$ and $\text{sat}_\Delta(x) = x$ for any $|x| \leq \Delta$; $e^L_{i,t}$ is the relative tracking error corresponding to the leader train 0 with $d_0$ standing for the minimal/safe gap; $f^{re1}_{i,t}$ denotes the basic rolling resistance per unit train weight; $\gamma_{i,t}$ is a given compensation function; $\Delta$ is the saturated control input bound; and $k_1, k_2, k_3$ are some positive gains to be determined.

In [58], the problem of velocity- and input-constrained tracking control of discrete-time high-speed train systems was solved. The following distributed cooperative tracking control law was proposed for each train using the local information from its adjacent trains:

$$u_{i,t} = \sum_{j\in\mathcal{V}} a_{ij}(x_{j,t} - x_{i,t} + d^{de}_{ij}) - k_{1,i}(v_{i,t} - v_{0,t}) + f_{i,t} \quad (4.5)$$

where $k_{1,i}$ is a constant damping gain. Then, based on the model transformation and convex analysis techniques, it is shown that all trains can eventually track the virtual leader train in a stable mode.

In [61], the following coordinated time-varying low gain feedback control law was proposed:

$$u_{i,t} = \rho B^T P(\theta_t) \sum_{j\in\mathcal{V}} a_{ij}(\tilde{s}_{j,t} - \tilde{s}_{i,t}) - a_{i0}\tilde{s}_{i,t} \quad (4.6)$$

for a multiple high-speed train system under unknown input delays:

$$\dot{\tilde{s}}_{i,t} = A\tilde{s}_{i,t} + B u_{i,t-\tau} \quad (4.7)$$

where $\tilde{s}_{i,t} = [x_{i,t} - x_{0,t} + (i-1)d_0, v_{i,t} - v_{0,t}, z_{i,t}]^T$; $z_{i,t} = F^{db}_{i,t}/((1+\gamma)M_i) - (k_{1,i} + k_{2,i}v_{i,t} + k_{3,i}k_{4i}v^2_{i,t})$; $\tau \leq \bar{\tau}$ is the unknown constant information transmission delay in the control channel; $A = [0, 1, 0; 0, 0, 1; 0, 0, 0]$, $B = [0; 0; 1]$; $\rho$ is a positive coupling gain; $P(\theta_t)$ is a positive definite matrix satisfying the following Algebraic Riccati:

$$A^T P(\theta_t) + P^T(\theta_t)A - P^T(\theta_t)BB^T P(\theta_t) = -\theta_t P(\theta_t), \forall \theta_t > 0. \quad (4.8)$$

Moreover, $\theta_t$ can be designed as $\theta_t = h/\hat{\tau}_t$ with $h$ being a small positive constant and $\hat{\tau}_t > 0$ being the estimated delay function for the unknown delay $\tau$. Then, it is shown that the coordinated control problem of VC-based trains is eventually converted into a stabilization problem for a linear multiple-input-multiple-output system subject to an unknown constant input delay.

Zhao and Wang [75] used the Artificial Potential Fields (APF) concept to design a controller for multi-train cooperative driving on urban railways:

$$u_{i,t} = -k_1 \sum_{j\in\mathcal{V}} a_{ij}\nabla \ln[\cosh(x_{j,t} - x_{i,t} - d^{de}_{ij})]$$
$$-k_2 \sum_{j\in\mathcal{V}} a_{ij}\nabla \ln(v_{j,t} - v_{i,t})$$
$$-k_3\nabla \ln[\cosh(v_{i,t} - v^{de}_i)] - k_4\nabla \ln[\cosh(v_{i,t})] \quad (4.9)$$

where $v^{de}_i$ the desired speed of the $i^{th}$ train; and $\nabla$ is the gradient of the potential. The controller includes four parts: the first part is for distance control; the second part is for speed consensus control; the third part is for target speed control; and the last part is for train braking control to reduce train speeds during brake applications. The first two parts involve T2T communication topologies in terms of $a_{ij}$.

The above reviewed controllers have all considered a flexible communication topology by including $a_{ij}$. A simple and special case of the communication topology is the Predecessor-Following topology where a follower train aims to follow a single leader right in front of it. For such a case, the controller can be significantly simplified. Liu [76] presented a simple train following controller as:

$$u_{F,t} = -k_1(d_{F,t} - d^{de}_{F,t}) - k_2[k_3 + k_4\tanh(k_5 d_{F,t-\tau} - k_6)] \quad (4.10)$$

where $u_{F,t}$ is the controller input for the follower; and $\tau$ is the communication delay. The controller considers two gains. The first one is for train position errors. The second one is for the following distance. Zhou et al. [77] also presented a train-following controller as:

$$u_{F,t} = -k_1\tanh(e^{sp}_{F,t}) - k_2\tanh(e^{sv}_{F,t}) - k_3 v_{F,t}/(v_{F,t} - v_{max}) \quad (4.11)$$

where $v_{max}$ is the speed limit. This controller incorporates both position and speed errors. However, they are inserted in a hyperbolic tangent function to smooth the control. An additional gain for the control of maximum speed is also included.

B. *Optimal Control (mainly Model Predictive Control (MPC))*

Optimal control is a group of control methods designed to find a control action for a dynamic system over a period (finite or infinite horizon) such that an objective is, or multiple objectives are optimized. Usually, an optimization process is involved to reach the optimal solution. For VC studies, the most common version is MPC (or moving finite horizon) [43,78] as shown in Fig. 7. Infinite horizon optimal control was used in [38] and a generalized predicative control (finite-horizon unconstrained) was used in [48,49].

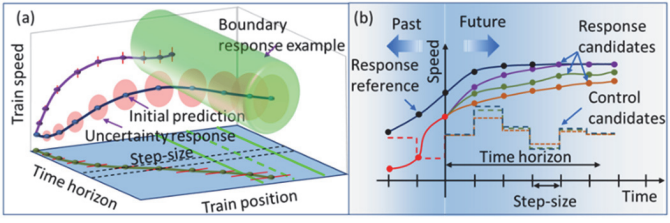

Fig. 7. Schematic diagram of model predictive control (MPC): (a) system prediction with uncertainty, and (b) selection from an optimization process

Various train motion models reviewed in Section 2 can be used to update each train's motion states. A fitness function or a cost function is then needed to assess the optimality of a specific solution. MPC for VC studies usually formulates the optimization problem as a minimization problem, i.e., the optimal solution has the minimum cost value. A basic cost function that considers one parameter can be expressed as:

$$J_{i,k} = \sum_{j=1}^{N_p} M_i \left( v_{i-1,k+j|k} - v_{i,k+j|k} \right)^2 \quad (4.12)$$

where $J_{i,k}$ is the cost function value for each train $i$ at current time step $k$; $N_p$ is the total number of prediction steps, i.e., the horizon length; and the subscript $j|k$ indicates $j$-th the predicted value at the $k$-th actual time step. This cost function was used by Xun et al. [65] to study the overspeed issue of virtually coupled trains. Su et al. [53] removed the mass term and replaced the train speeds by train positions to conduct VC simulations. In [54], three cost functions were used for acceleration, speed, and position. They all had the same format as the one used in [53]. The position cost function was defined as a 1-norm function whilst the other two were 2-norm functions. Three cost functions were not combined but used individually according to different train operational conditions. Both [53] and [54] used a Neural Network model to predict the dynamics of a preceding train.

When multiple parameters are considered, the cost function can be written as:

$$J_{i,k} = \sum_{j=1}^{N_p} \left[ k_1 \left( e_{i,k+j|k}^{sp} \right)^2 + k_2 \left( u_{i,k+j|k} \right)^2 + k_3 \left( d_{i,k+j|k} \right)^2 \right] \quad (4.13)$$

$$J_{i,k} = \sum_{j=1}^{N_p} \left[ k_1 \left( e_{i,k+j|k}^{sp} \right)^2 + k_4 \left( \Delta u_{i,k+j|k} \right)^2 \right] + k_5 \left( e_{i,k+N_p|k}^{sp} \right)^2 \quad (4.14)$$

where $k_1 - k_5$ are weighting factors. The three elements of the first equation correspond to gap errors, actuator forces/powers, and train following gaps, respectively. In the second function, the first element is the gap error whilst the second and third elements are actuator force/power variations and the gap error of the very last prediction point. This literature survey shows that all MPC cost functions for VC simulations consider the gap error element. This is understandable as minimum gap error is the primary objective of the VC operation. Luo et al. [34] added the actuator force/power elements. This represents an energy saving objective for the actuator; it also helps to keep controller forces within the system limits. In [34], an estimator was used to update control law and to minimize the differences between predicted train states and actual train states. Researchers in [44,49,57] used the gap error element and the actuator force/power variation element without the force/power element. The variation term not only served a similar role as the force/power element but also helped to minimize the jerk, i.e., the variations of forces/accelerations. Wu et al. [50] directly used a jerk (variation rate of acceleration) element in the cost function. She et al. [56] used the gap error element, force variation term, and following gap element. The utilization of the gap element reflected the objective of minimum headway and could help to increase the traffic efficiency. Liu et al. [52,55] used the gap element, force/power element and the final prediction element. The final prediction element could be regarded as a terminal term which reflected the control quality beyond the predicted horizon. To some extent, the use of the final prediction element extended the prediction to an infinite horizon. Liu et al. [52,55] also considered speed differences between two compared trains; the speed difference can be used as the gap element. Felez et al. [43] used a different cost function which also considered following gap, but added speed differences and jerk limits:

$$J_{i,k} = \sum_{j=1}^{N_p} \left( k_1 \left| 1 - d_{i,k+j|k} / d_{i,k+j|k}^{de} \right| + k_6 \left| e_{i,k+j|k}^{sv} / e_{max}^{sv} \right| + k_7 \left| J_{i,k+j|k}^{f} / J_{max}^{f} \right| \right) \quad (4.15)$$

where $k_6$ and $k_7$ are also weighting factors; $e_{max}^{sv}$ is the speed difference limit; $J_i^f$ is the jerk; and $J_{max}^f$ is the jerk limit.

Liu et al. [38] formulated MPC as an infinite horizon control problem in which the cost function was expressed as a continuous function:

$$J_{i,t} = \int_0^\infty exp(-k_8 t) \left[ k_1 \left( e_{i,t}^{sp} \right)^2 + k_2 \left( e_{i,t}^{sv} \right)^2 + k_3 \left( e_{i,t}^{v} \right)^2 + 0.5 \left( u_{i,t} \right)^2 \right] dt \quad (4.16)$$

where $t$ is time in the predicted horizon; $k_1 - k_3$ are weighting factors; $e_i^v$ is the speed difference between the assessed follower with the platoon leader; and $k_8$ is a discount factor to adjust the terminal cost. Knowing that the $exp(-k_8 \tau)$ function has a falling effect with the increases of time, predicted train performance far into the future contributes less than that in the near future. Su et al. [51] also used the exponential function, and the cost function was expressed as:

$$J_{i,t} = 0.5 \int_0^{T_p} exp(-k_8 t) \left[ k_1 \left( e_{i,t}^{sp} \right)^2 + k_2 \left( u_{i,t} \right)^2 \right] dt + 0.5 exp(-k_8 t) k_1 \left( e_{i,T_p}^{sp} \right)^2 \quad (4.17)$$

where $T_p$ is the time length of the predicted horizon. This function has considered the terminal cost by adding $e_{i,T_p}^{sp}$.

It is noted that all the cost functions above have a decentralized format, which means that the assessed local platoon only includes the assessed train itself and its leader(s). All the relative error terms are to be calculated between the assessed train and its leader(s). For most cases, there is one leader, and the leader is the adjacent front train. For (4.16), speed difference between the platoon leader $e_{i,t}^v$ was also considered. If the cost functions are written in a centralized format, the cost function values are to be summed across all assessed trains. Taking (4.12) as an example, its centralized format can be expressed as:

$$J_k = \sum_{j=1}^{N_p} \sum_{i=2}^{N} M_i\left(v_{i-1,k+j|k} - v_{i,k+j|k}\right)^2. \quad (4.18)$$

The literature review shows that most of the MPC controllers for VC studies were decentralized. Centralized MPC control can be found in [51,65]. An advantage of using the centralized format was that it allows better analyses of global string stability [51].

The computing time and the length of the prediction horizon are important issues for MPC. On the one hand, to be implemented in real-world operations, the VC controller needs to be faster than real-time and must allow some time margins for the executions of other necessary steps, such as communications and actuator executions. The overall control process must be real-time. On the other hand, MPC prefers a longer prediction horizon, i.e., more prediction steps to increase the performance of the controller. Longer prediction horizons commonly mean longer computing time. Researchers in [53, 54, 63] reported increases in both controller performance (less tracking errors) and computing time (see Fig. 8) with the increases of prediction steps. Most of the MPC applications for VC studies used no more than 10 prediction steps: Su *et al.* [51] used 10 prediction steps; researchers in [52,54,55,57,65] used 5 prediction steps; researchers in [34,50,56] used 4 and 3 prediction steps respectively.

### C. Sliding Mode Control (SMC)

SMC is a nonlinear discontinuous controller, as shown in Fig. 9, which controls the dynamics of a system by using a discontinuous control signal. The direct control objective is to drive the system to slide along a surface that satisfies certain system state requirements (speed, position, acceleration, etc.) and eventually reach an equilibrium. There are two major components in SMC: the sliding surface (manifold) and the controller.

In the case of VC control, a simple sliding surface can be designed by considering the position errors and speed errors:

$$S_i = e_{i,t}^{sp} + \kappa_i\, e_{i,t}^{sv}, \quad \forall\, i \in \mathcal{V} \quad (4.19)$$

where $\kappa_i > 0$ is a constant to be chosen; and $S_i$ is the sliding variable. The sliding variable can be interpreted as a compound error parameter of the system, which also considers a weighting factor ($\kappa_i$) between different error components. If needed, an acceleration error can also be added, however, the VC operation is more focused on speed and position errors. If the desired control input $u_{i,t}$ exists for a system such that $\lim_{t \to \infty}|S_i| = 0$, then the regulation of the spacing errors and speed errors can be achieved, namely, $\lim_{t \to \infty}|e_{i,t}^{sp}| = 0$, and $\lim_{t \to \infty}|e_{i,t}^{sv}| = 0$.

As shown in Fig. 9, SMC usually has a reaching phase. Liu *et al.* [45] added a nonlinear function to the sliding surface which achieved better robustness by avoiding the reaching phase and making sure the initial state of the system is on the sliding surface

$$S_i = \left(e_{i,t}^{sp} - r_t\right) + \kappa_i\left(e_{i,t}^{sv} - \dot{r}_t\right), \quad \forall\, i \in \mathcal{V} \quad (4.20)$$

where $r_t$ is a 5th order piece-wise polynomial related to initial position error, velocity error, and acceleration error. The controller was eventually expressed as:

$$u_{i,t} = -\frac{1+\gamma}{0.001g}\left[k_1\left(e_{i,t}^{sv} - \dot{r}_t\right) + f_{i,t} - \ddot{x}_{0,t} - \ddot{r}_t + k_2 \mathrm{sat}(S_i)\right] \quad (4.21)$$

where $k_1$ and $k_2$ are constant gains. To handle the chattering issue, the control action was smoothed as:

$$\mathrm{sat}(S_i) = \begin{cases} 1 & S_i > \Delta \\ \frac{S_i}{\Delta} & |S_i| \le \Delta \\ -1 & S_i < -\Delta \end{cases} \quad or \quad \mathrm{sat}(S_i) = \tanh(\Delta S_i) \quad (4.22)$$

where $\Delta$ is a smooth control parameter.

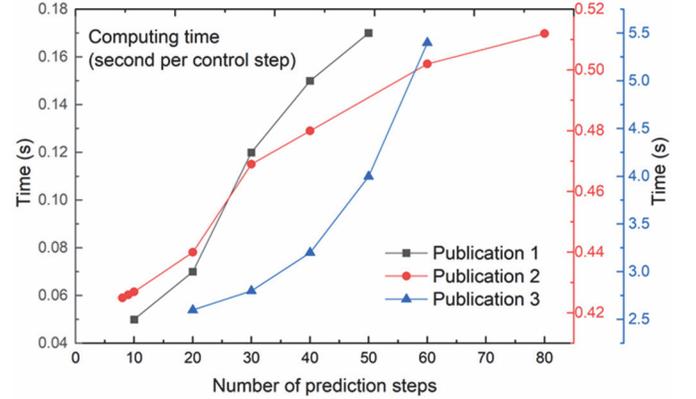

Fig. 8. Influences of prediction steps on computing time of MPC (Publications 1-3 are [53,43,44])

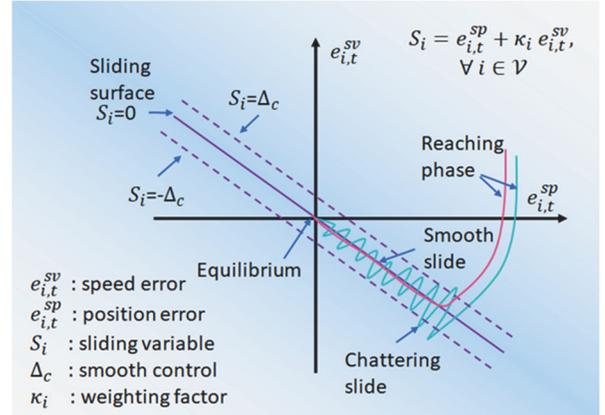

Fig. 9. Schematic diagram of sliding mode control (SMC)

Another SMC application can be found in [47] where two functions called position APF and speed barrier function were used to design a sliding surface as:

$$S_i = (k_1 + L_i)\tanh\left(e_{i,t}^v\right) + (k_2 + W_i)e_{i,t}^p \quad (4.23)$$

$$W_i = \min\left\{\left(k_{3,i}\frac{|d_{i,t}|^2 - d_{i,t}^{de^2}}{|d_{i,t}|^2 - d_{i,t}^{min^2}}\right), 0\right\}^2, \quad L_i = \frac{k_{4,i}^{\,2}}{\pi}\tan\left(\frac{\pi\left(e_{i,t}^v\right)^2}{2k_{4,i}^{\,2}}\right) \quad (4.24)$$

where $k_1$-$k_4$ are positive parameters to be adjusted; $L_i$ is the speed barrier function; $W_i$ is the position APF; and $d_{i,t}^{min}$ is the minimum gap. This sliding surface is also fundamentally based on speed and position errors of the follower train. However, it incorporates a position constraint (minimum gap) and a speed constraint (the leader speed). The study [47] also considered uncertain parameters, disturbances, and control input saturations (constant limits). The controller is regarded as adaptive and stable if external forces and uncertainties are within the bounded values.

In [46], a robust SMC gap controller was proposed for two virtually coupled railway trains. Vehicle mass errors, acceleration limits; velocity dependent disturbances; and time-variable disturbances were bounded and processed. The system is regarded as stable if the disturbances and uncertainties are within the bounded values. To reduce disturbances caused by incorrect information about wheel diameters to VC controllers, Park et al. [46] proposed to use balise information to estimate and correct speed and position errors. Simulations show that, by using the error corrections, the gap errors during VC operations are reduced from mostly 5 m to less than 1 m.

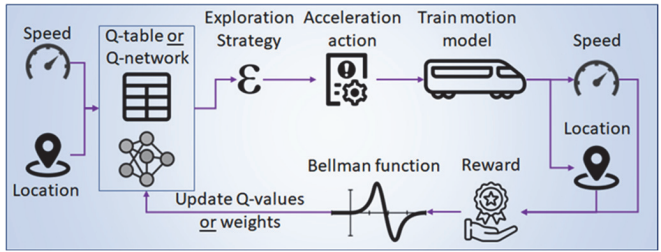

Fig. 10. Training process for VC controller based on Q-learning and Deep-Q-Network

### D. Machine Learning (ML)-Based Control

ML is an emerging research field which includes several powerful techniques. For VC studies, the Q-learning [39] and Deep-Q-Networks (DQN) [37] techniques have been used. The Q-learning technique is a model-free ML technique which learns the value of an action (T/B forces) in a particular state (speed and position errors). It is based on the finite Markov Decision Process; and the objective is to find a policy that can maximize the reward value over any/all successive steps. In [39], a Q-learning-based cooperative VC approach was developed as shown in Fig. 10. To preserve the global optimality, an APF was used to measure the effects of actions and derive the reward function:

$$rwd_i = -k_1 \sum_{j=1}^{N} a_{ij} \nabla \ln\left(\cosh\left(e_{i,j}^{sp}\right)\right) - k_2 \sum_{j=1}^{N} a_{ij} \tanh\left(e_{i,j}^{sv}\right) \quad (4.25)$$

where $rwd_i$ is the reward value; $k_1$ and $k_2$ are two positive constants to determine the weights of two errors; and $\nabla$ is a gradient of the potential. This reward function considered the platoon communication topology via $a_{ij}$ as discussed in Section 4.1. The Q-learning method was selected to obtain the cooperative control policy for each train in the platoon. This means the values that got updated were Q-values in Fig. 10 and the Q-value structure used was the Q-table. In this study, acceleration actions were limited by maximum and minimum values. Twenty-one and 22 states were used for position error and speed error respectively in the Q-table.

In [37], a DQN algorithm was used to solve cooperative collision-avoidance control for two virtually coupled trains. The application of DQN was similar to that of Q-learning as shown in Fig. 10. In the case of DQN, the values that got updated were the weights of the Neural Network and the Q values were stored in the structure of a Neural Network rather than a Q-table. This DQN algorithm integrated the artificial Neural Network with the reinforcement learning technique; the combination could overcome the Curse of Dimensionality in the Q-learning. In this study, it was assumed that the operation trajectory of the preceding train was predictable and could be obtained by the following train in real time through T2T communications. The reward function was expressed as:

$$rwd_{i,k} = \begin{cases} k_1 & v_{i,k+1} > v_{max} \text{ or } d_{i,k+1} < D_{safe} \\ d_{i,k} - d_{i,k+1} & 0 \le v_{i,k+1} \le v_{max} \text{ and } d_{i,k+1} \ge D_{safe} \\ k_2 & D_{safe} < d_{i,t}^{min} \text{ and } d_{i,k+1} < d_{i,t}^{min} \end{cases} \quad (4.26)$$

where $k_1$ and $k_2$ are two negative rewards (punishments); $k$ indicates the current control step; and $d_{i,t}^{min}$ is the minimum headway. This reward function was constrained by safety limits that made sure the speed did not exceed maximum speeds to avoid collisions. This study used speed dependent maximum T/B forces. The Neural Network had 4 hidden layers and a maximum 30 nodes in each layer.

Two other applications that were related to ML based techniques can be found in [53, 54]. These two applications were not based on ML as the VC simulations were achieved by using MPC as reviewed in Section 4.2. The long and short-term memory (LSTM) Neural Network was used to predict the dynamics of the preceding train.

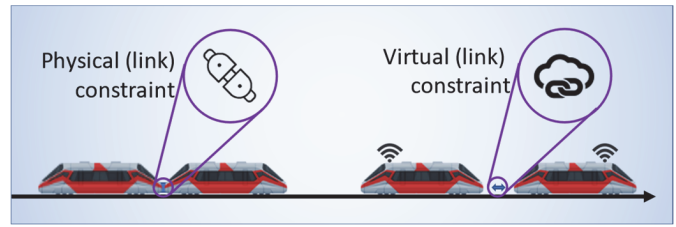

Fig. 11. Physical constraint and virtual constraint between two coupled trains

### E. Constraints-Following Control (CFC)

Traditionally, when two trains need to operate together as a train set, a physical link is used to couple these two trains as shown in Fig. 11. The physical link, from the system dynamics point of view, acts as a physical constraint. The constraint forces two coupled trains to have (as much as possible) constant distance, zero relative speed and zero relative acceleration. Adopting the constraint idea to VC controller designs, one can design a controller to act as a virtual constraint between two virtually coupled trains. In the VC case, the constraint forces that are needed to make the follower train comply with distance, velocity and acceleration constraints are then the controller forces needed to be applied on the follower train. This is the basic idea of constraints-following control. Zhang et al. [62] and Wang et al. [63] and used this type of controller for VC studies. The acceleration constraint was expressed as (3.3) while the speed and position difference constraints were expressed as:

$$e_{i,t}^v = k_3 \exp(-k_4 t) + e_{i,t}^{vde}, \quad d_{i,t} = k_5 \exp(-k_6 t) + d_{i,t}^{de} \quad (4.27)$$

where $k_3$-$k_6$ are constants to adjust the transition process; and $e_{i,t}^{vde}$ is the final desired speed difference.

To implement these constraints, an unconstrained dynamic system was first modelled. Then the constraints were imposed to the system model by following the principle of virtual displacement and the d'Alembert principle. Eventually, the controller was expressed as:

$$u = M^{1/2}(HM^{1/2}) + [b + HM^{-1}(C\dot{x} + G)] \quad (4.28)$$

where $u$ is the constraint force, i.e., controller input; $M$ is the mass matrix; $H = [1, -1]$; $b$ is the acceleration constraint; $x$ is the train position; $C$ consists of all speed-related terms in the rolling resistance formula; and $G$ represents grade forces and rolling resistance that are independent from train speeds. Both controllers [62,63] did not apply limitations to maximum T/B forces. Zhang et al. [62] added the consideration of uncertainties for train masses and external forces.

## V. COMPARISONS AND FURTHER DISCUSSIONS OF VC CONTROL TECHNIQUES

In this section, the advantages and disadvantages of the existing controllers for VC applications, including CBC, MPC, SMC, ML, and CFC, are discussed.

### A. Consensus-Based Control (CBC)

The primary advantage of CBC is that it explicitly considers T2T communication topologies used in VC applications. This enables better and easier analyses for various VC topics by using the tool of Graph Theory, optimal dynamic T2T communication topology and communication failure tolerance as well as network-induced constraints (including communication delays, data losses, fading channels, varying topologies). Another advantage of CBC is that it can be very simple. For example, when the predecessor-following topology is used, the controller becomes a simple feedback controller. Simplicity is beneficial for quicker uptakes of VC studies and applications.

Having understood the CBC's advantages in diverse communication topologies, it is also argued that VC applications often have much fewer agents (trains) than many other multi-agent systems such as a swarm of drones and mobile robots. Currently, three or four trains may be considered as the maximum number of trains in a practical platoon. Under this situation, complex communication topologies may be unnecessary in VC applications. However, communication topology analyses will be interesting when smaller train units are discussed, such as, passenger pods for Hyperloop [79] and ultra-light people movers [80]. These train units have the size of passenger cars and may need to dispatch in large numbers during rush hours.

An obvious disadvantage of CBC is that it has limited capability to handle constraints. This is also reflected by the review in Section 3.2. Most of the CBC studies did not consider any or consider only a few constraints. It is noted that some constraints, such as the saturation constraint (maximum forces/accelerations), are important parts of realistic systems. Constraints that were not designed in controllers but were encountered during operations may cause system disabilities.

### B. Model Predictive Control (MPC))

The prominent advantage of MPC is probably its great capability of constraint handling. Various constraints can be applied by modifying the final values of the cost functions via checking the candidate results against the constraints. Soft constraints can also be integrated into the cost function itself such as the jerk limit in (4.15). The review in Section 3.2 shows that applications where MPC was used were able to consider almost all constraints that were reviewed. These constraints can help to develop realistic system models and increase the applicability of the final controller. Also due to the utilization of cost functions, multi-objective optimizations are possible with MPC. Controller inputs can have good performance in terms of multiple aspects such as ride comfort and energy consumptions.

The prediction feature of MPC is also beneficial for train driving as railway trains usually have heavy masses; driving commands need to not only consider the current train status but also to anticipate a certain length of track section in advance of the current train position. This is especially the case for freight trains. Many Automatic Train Operation (ATO) functions for freight trains have the so-called 'look-ahead' function [81] by using which a prediction function was achieved and considered while determining the current train driving commands.

MPC also has an obvious disadvantage which is its high computational cost. MPC requires multi-step predictions and an optimization process for each prediction. The prediction feature means that, for example, if the number of prediction steps is five, the computing time of MPC is most likely five times slower than controllers that do not need predictions. In addition, the optimization process is also time consuming, especially when multiple constraints and multiple objectives are involved. More constraints mean longer computing time; nonlinear constraints such as speed-dependent maximum forces and robust MPC [38,78] take even longer computing time. Results from Figure 8 have shown that some MPC applications were not able to achieve real-time inputs. Simulation cases in [52] indicated that, when the controller input frequency is higher than five times per second, MPC is difficult to achieve real-time inputs.

### C. Sliding Mode Control (SMC)

SMC is a robust discontinuous controller that has a proven equilibrium point which can be reached in finite time when unconstrained. SMC also has good robustness for complex systems with high-order nonlinearities. For VC applications, SMC has been described as having small tracking errors and good performance while handling uncertainties. System uncertainties are a significant issue in VC applications. The uncertainties can originate from almost all force elements of train motion models. For example, train masses can vary between rush hours and non-rush hours due to different numbers of passengers onboard; resistance forces are well known for their uncertainties when estimated by using empirical formulas [30]; and T/B forces also have uncertainties resulting from friction phenomenon. Therefore, the good uncertainty handling capability of SMC is certainly an advantage, which has been well demonstrated by [46, 47, 82].

It is also noted that, while deriving SMC's proven convergence to the equilibrium point, it was assumed that the control cost was not limited. This can be unrealistic as actuator forces, i.e., T/B forces, always have limitations. For VC applications, T/B forces may also need to be adjusted to achieve better ride comfort and/or lower energy consumptions. When force limitations and various other constraints are applied, there are possibilities for the system to become unstable. In addition, the robustness property of SMC also becomes unprovable when chattering mitigations are applied, see (4.22). Without the chattering mitigations, the system will be constantly fluctuating around the equilibrium and railway train actuators are not able to execute the high frequency controls delivered from SMC.

### D. Machine Learning (ML)-Based Control

ML-based VC controllers naturally inherit the 'model-free' advantage of all data-driven approaches; detailed train system models are not needed anymore. Controller actions are derived from an Artificial Intelligence structure (e.g., Q-network or Q-table) which are developed from the training materials, i.e., historical data.

Due to the model-free feature, the final form of ML-based controllers is simple, which offers the advantages of simple computer implementations and fast computing speeds. These are important as onboard computing resources are very limited and real-time control is essential for real-world implementations. For example, Basile et al. [24] reported that ML-based controllers could achieve similar performance as MPC controllers for road vehicle platoons. In this case, the former outperformed the latter in terms of computing speed, especially in more realistic and complex scenarios with multiple nonlinear constraints. It is noted that multi-objectives and constraints can be incorporated into the reward function of the training process.

According to [24], ML-based VC controllers that had an accurate learning process should be able to reach the optimal solution in real-time considering all current uncertainties and contingencies without any prior knowledge of the railway environment. In other words, ML-based VC controllers should have advantages in dealing with system uncertainties.

However, it is also argued that the performance of ML-based controllers heavily depends on the training process. These methods work in high-dimensional action spaces that are difficult to explore efficiently; any datasets would be inadequate for this purpose [24]. There are also unforeseen situations which cannot be fully covered by the training materials. In addition, VC is still an emerging technology, training data that can be used for controller development is very limited. These factors could inevitably influence the performance of ML-based VC controllers.

Finally, ML-based VC controllers also share the same concern of other data-driven approaches for using black-box models. The controller structure does not reflect the actual system; the stability and performance of the system cannot be theoretically examined. There are understandable concerns for ML-based models to be used for safety-critical controls such as the VC operations.

### E. Constraints-Following Control (CFC)

CFC is a control method based on the Udwadia-Kalaba formulation which was designed to derive s of Motion (EoM) for constrained mechanical systems. In this sense, CFC has a clear physical meaning of replacing physical constraints with virtual constraints as shown in Figure 11. It is noted that these constraints are different from the system constraints discussed in Section 3.2. The former is more of a connection or link between two trains whilst the latter is more of system limitations. Theoretically, unconstrained (unlimited) system stability of CFC can be proved by following Lyapunov stability theory. However, it is not guaranteed when system constraints (limitations) such as maximum T/B forces are applied. Neither [62] nor [63] considers any system constraints (limitations). It is noted the minimum following distance can be considered by adding inequality components to the virtual constraints and maximum relative motions (distance, speed, and acceleration) between the follower and the leader can be controlled by designing the expressions of the virtual constraints. Zhang et al. [62] also demonstrated the consideration of bounded uncertainties into the controller.

## VI. FUTURE RESEARCH TOPICS

Future research topics are discussed from two related perspectives. The first one is for better controller development whilst the second one is for better implementation of VC controllers.

### A. For Better Controller Development

There is a list of topics that can be further studied to help to develop more robust, more realistic, and more accurate VC controllers.

1) As reviewed and discussed in Section 2.2, rotational inertia has considerable influences on train motion simulations. However, very few publications have considered rotational inertia in train motion models. This point can be improved for future controller development.
2) As reviewed and discussed in Section 3.1.1 (see Figure 4), a transition is desired for any switches between traction and brake to reduce longitudinal impacts. However, only one publication [46] was found that considers such transition. Therefore, future work can be conducted towards this direction to design optimal transition profiles for VC merging and splitting.
3) Energy consumption as an important topic has been widely addressed in ATO related studies [83]. However, current VC operations often design follower trains to have the same speed as the leader trains. Due to location differences, leader speed profiles may not be optimal anymore for follower trains. Therefore, energy optimal VC operation is a research topic that of great research and engineering values.
4) When multiple objectives are considered for cost functions of MPC and reward functions of ML-based controllers, different expressions and different weighting factors can be explored. To the best of the authors' knowledge, such studies have not been reported for VC operations.

5) As reviewed in Section 3.2, most maximum force/acceleration constraints use constant limits. However, maximum railway train T/B forces are strongly speed dependent. Often, maximum T/B forces at higher speeds are significantly lower than those at lower speeds. Further research into achieving speed dependent maximum T/B force constraints is an important topic.

6) Communication safety is critical as VC operations heavily rely on wireless communications [84]. Controllers that are safe and robust against communication faults and communication attacks are rare [85,86]. Further research along this direction is necessary for successful implementations of VC technologies.

7) VC control and operation under different communication topologies and/or optimal communication topology have not been adequately studied. It is also promising to consider a co-design of VC controller and T2T communication topology for realistic VC applications. On the other hand, most VC controllers are designed for only single train platooning, while VC controllers for multiple train platoons have not been available.

8) Current VC studies still have a focus on homogeneous train platooning. Heterogenous VC studies where different trains have different masses, traction performance, and braking performance [87] are interesting and of good research and engineering values.

9) Studies about VC controller robustness against actuator accuracies are not seen in the open literature. Railway train T/B forces need to be eventually generated from wheel-rail contact. There are many variable steps between T/B commands generated from VC controllers and actual T/B forces generated from wheel-rail contact. This can be an important topic towards safe and successful implementation of VC operations.

10) Intermittent or low frequency train controllers which require less frequent actuator actions, e.g., once every two seconds, are more suitable for large systems such as railway trains. Intermittent controllers also agree with traditional manual train driving. However, intermittent, or low frequency train controllers for VC operations have not been published.

11) Dynamic Coupling has been proposed [88,89] to allow virtually coupled railway trains to perform physical coupling and decoupling at cruising speeds. Dynamic Coupling controllers can be based on VC controllers but have much higher requirements in terms of train motion controls. This is still a new topic.

12) VC operations across junctions (switches or turnouts) have been identified as slow points for traffic efficiency improvement [90,91]. A possible solution is to design a controller to form VC before the junction and split after the junction. Such an issue deserves deep investigation.

13) Multiple-carriage train models and in-train forces have been considered in position and velocity tracking control problems for individual trains [82]. How to consider in-train forces in VC operations is still challenging and an interesting research topic.

B. *For Better Controller Implementation*

The following topics can be further studied to reduce the requirements for VC controllers and help to achieve easier and better implementations of VC controllers.

1) Train operational resistances including rolling resistance and curving resistance vary considerably depending on train speeds, train conditions, track conditions and environment conditions. Empirical formulas used for VC controllers have been known as inaccurate. Online resistance force assessments can help to deliver better resistance inputs for controllers and achieve more accurate VC controls.

2) More accurate resistance formulas can be developed. This can be progressed by conducting field measurements and analyzing historical data. Detailed high-fidelity models can also help to develop more accurate formulas.

3) Actual T/B forces generated from wheel-rail contact can be considerably different from controller commands due to various factors including wheel-rail adhesion uncertainties, brake friction variations, and system response delays. T/B force determination and feedback can help to achieve better VC controls. For example, instrumented couplers can be used to measure the actual traction or dynamic brake forces delivered from motor vehicles to trailing vehicles.

4) Train masses also represent another significant uncertainty for VC controllers. For example, a fully loaded passenger vehicle can weigh up to 50 tonnes whilst an empty vehicle weighs only 30 tonnes. Modern weigh-in-motion technology installed on railway tracks can help to reduce this uncertainty. Research can also be carried out to determine train masses by using suspension sensors or longitudinal acceleration sensors.

5) Track data is an important input for VC controllers and is needed while determining track gradient forces, curving resistance, and tunnel resistance. Accurate track data surveys or on-board information correction algorithms can help to achieve more accurate VC control.

6) Current VC studies differentiate T/B forces and traction and brake forces. However, brake forces can be generated by frictional brakes and dynamic brakes. Therefore, a further brake blending algorithm is needed to achieve realistic train controls under VC operations.

7) Real-time controllers are essential for real-world VC operations. When controller computing speed is low, parallel computing can be used to improve computing speeds. For example, when iterative optimization is required for MPC, parallel computing can be used to process multiple candidates concurrently so as to improve computing speeds.

8) As mentioned by Basile *et al.* [24], VC simulators can be developed to generate data to train ML-based models. The simulator can also be used to validate and verify various parts of VC operations such as controller performance and train dynamics.

## VII. Conclusions

VC is an advanced technology that has great potential to significantly improve railway traffic efficiency. VC control techniques are critical developments towards successful implementations of this technology. This paper provided a systematic review regarding the existing and emerging control techniques used for VC studies. Train motion models that are necessary for model-based controller development are reviewed first. Force elements considered in the train motion models are then examined and classified into T/B forces, rolling resistance, tunnel resistance, curving resistance and grade forces. Next, the relative train motion control objectives and constraints for VC controllers are discussed in detail. Various minimum gap reference expressions are provided as well as different transitional gap references. Furthermore, constraints used for VC controllers are presented, which includes minimum gap, maximum speed, maximum T/B forces (accelerations), maximum powers, and maximum jerks. An emphasis is then placed on five emerging controllers used for VC studies, which are classified into CBC, Optimal Control/MPC, SMC, ML-based control, and CFC. The advantages and disadvantages of these controllers for VC applications are also elaborated. Finally, two lists of topics for future studies regarding controller development and controller implementation respectively are introduced.

## Acknowledgments

The editing contribution of Mr. Tim McSweeney (Adjunct Research Fellow, Centre for Railway Engineering) is gratefully acknowledged.